\documentclass{article}
\usepackage{epsfig}
\usepackage{amssymb}
\def\ds{\displaystyle}

\newtheorem{theorem}{Theorem}[section]
\newtheorem{coroll}{Corollary}[section]
\newtheorem{lemma}{Lemma}[section]
\newtheorem{remark}{Remark}[section]
\newtheorem{proposition}{Proposition}[section]
\newtheorem{definition}{Definition}[section]
\def\le{\left}
\def\m{\mathop}
\def\ri{\right}
\def\br{\begin{remark}\rm\small}
\def\1{{\bf 1}}
\def\er{\end{remark}}
\def\bt{\begin{theorem}\rm}
\def\et{\end{theorem}}
\def\bc{\begin{coroll}\rm}
\def\ec{\end{coroll}}
\def\bl{\begin{lemma}\small}
\def\el{\end{lemma}}
\def\bd{\begin{definition}}
\def\ed{\end{definition}}
\def\bp{\begin{proposition}\rm}
\def\ep{\end{proposition}}
\def\be{\begin{equation}}
\def\ee{\end{equation}}

\def\bea{\begin{eqnarray}}
\def\eea{\end{eqnarray}}
\def\beas{\begin{eqnarray*}}
\def\eeas{\end{eqnarray*}}
\def \pa{\partial}
\def\C{{\mathbb C}}
\def\L{{\mathcal L}}
\def\R{{\mathbb R}}
\def\N{{\mathbb N}}
\def\Z{{\mathbb Z}}
\def\a{{\alpha}}
\def\b{{\beta}}
\begin{document}

\begin{flushright}
CRM-2842 (2002)
\end{flushright}
\begin{center}
\bf \Large 
\fontfamily{cmss}
\fontsize{17pt}{27pt}
\selectfont
 Bilinear semi-classical moment functionals and their integral
representation\footnote{ Work supported in part by the Natural
Sciences and Engineering Research Council 
of Canada (NSERC) and the Fonds FCAR du Qu\'ebec.}
\end{center}
\vspace{0.5cm}
\begin{center}
\large Marco Bertola$^{\dagger\ddagger}$\footnote{e-mail: bertola@crm.umontreal.ca}
\end{center}
\smallskip
\begin{center} 
$^\dagger$ Centre de recherches math\'ematiques,
Universit\'e de Montr\'eal\\ C.~P.~6128, succ. centre ville, Montr\'eal,
Qu\'ebec, Canada H3C 3J7\\
\smallskip 
$^\ddagger$ Department of Mathematics and
Statistics, Concordia University\\ 7141 Sherbrooke W., Montr\'eal, Qu\'ebec,
Canada H4B 1R6\\
\bigskip
\hrule
\bigskip
{\bf Abstract}
\end{center}
We introduce the notion of bilinear moment functional and study their
general properties. The analogue of Favard's theorem for moment
functionals is proven. The notion of semi-classical bilinear
functionals is introduced as a generalization of the corresponding
notion for moment functionals and motivated by the applications to
multi-matrix random models. Integral
representations of such functionals are derived and shown to be
linearly independent.\\
\\
Keywords:\ moment functionals, biorthogonal polynomials, semiclassical
functionals \\
\hrule
\section{Introduction}
The notion of moment functional is most commonly encountered as a
generalization of the context of Orthogonal Polynomials (OP) \cite{Szego}.
These are generally defined as a  graded polynomial orthonormal basis
in $L^2(\R,{\rm 
d}\mu)$ where ${\rm d}\mu$ is a given positive measure for which all
{\em moments} 
\be
\mu_i:= \int_\R{\rm d}\mu(x)\, x^i\ ,
\ee
exist finite. The moment
functional associated to such a measure is then the element  $\mathcal
L$ in the dual  space of polynomials, $\C[x]^\vee$ defined by 
\be
\mathcal L(p(x)) := \int_{\R}{\rm d}\mu \, p(x)\ ,\label{onemoment}
\ee
and it is uniquely characterized by its moments. The positivity of
the measure implies that we can always find orthogonal polynomials
which are real, so that the orthogonality relation reads
\bea
\mathcal L(p_m(x)p_n(x)) = h_n\delta_{nm}\ .\label{OPs}\\
p_n(x) = x^n+\mathcal O(x^{n-1})\,\in \R[x]\ ,\ \ h_n\in \R_+^\times.
\eea
Generalizing this picture one is led to consider {\em complex}
funtionals \cite{chihara}, i.e.  whose moments are not
necessarily real. The associated OPs are then defined by the same
relations (\ref{OPs}) where now the polynomials belong to the  ring $\C[x]$ and
$h_n$ are nonzero complex numbers.\par
One of the main applications of OPs is in the context of random
matrices \cite{mehta,RM2} where they allow to write explicit 
 expressions for the  correlation functions
of eigenvalues and  of the partition
function of these models.\par
Recently \cite{RM3,RM5,BEH,RM4} growing attention is devoted to the
$2$-matrix models (or the multi matrix models) in which the
probability space is the space of couples (or $n$-tuples) of matrices.
Also such models can be ``solved'' along lines similar to the one
matrix models by  finding certain bi-orthogonal
polynomials (BOP). The probability measure is given by 
\be
{\rm d}\mu(M_1,M_2) = \frac 1 {\mathcal Z_n} {\rm e}^{{\rm
Tr}(M_1M_2)}\,{\rm d}\mu_1(M_1)\,{\rm d}\mu_2(M_2)
\ee
where $M_i$ are $N\times N$ Hermitian matrices (usually) and the
positive 
measures  ${\rm d}\mu_i$ are $U(N)$ invariant.
The relevant BOPs are then a pair of graded polynomial bases
$\{p_n(x)\},\ \{s_n(y)\}$ ``dual'' to each other in the sense that 
\bea
\int_\R\int_\R{\rm d}\mu_1(x){\rm d}\mu_2(y)\, p_n(x)s_m(y){\rm e}^{xy} =
h_n\delta_{nm} \ ,\ \label{BOPs}\\
p_n\in \R[x],\ \ s_n\in \R[y]\ ,\ \  h_n\in \R^\times. 
\eea
The integral in Eq. (\ref{BOPs})  defines a particular kind of {\em bi}-moment
functional, that is an element of the dual to the tensor of two spaces of polynomials
$\C[x]\otimes_\C\C[y]$
\be
\mathcal L\big(p(x)\big|s(y)\big):= \int_\R\int_\R{\rm d}\mu_1(x){\rm
d}\mu_2(y)\, p(x)s(y){\rm e}^{xy}  \ ,
\ee
provided all its {\em bi-moments} $\mu_{ij}$ are finite 
\be
\mu_{ij}:= \mathcal L(x^i|y^j) \in \R\ .
\ee
Generalizing this picture we now consider complex bi-moment
functionals which are uniquely characterized by their (complex)
bi-moments $\mu_{ij}\in \C$.\par
The notion of semiclassical moment functional for a functional of the
form (\ref{onemoment}) requires that the measure ${\rm d}\mu(x)$ has a
density $W(x)$ whose logaritmic derivative is a rational function of
$x$ and the support is a finite union of intervals. This condition can
be translated into a distributional equation for the moment functional
itself and then generalized to the complex case
\cite{laguerre,maroni,shohat}.\\
Motivated by the applications to 2-matrix models, 
we are interested in the corresponding notion of semiclassical
bi-moment functionals (which we will define properly later on) and in 
studying their properties: we will produce (complex
path) integral representations for them, generalizing the framework of 
\cite{ismail,marce,marce1} to this situation.\par
We quickly recall that  \cite{laguerre,maroni,shohat} a  moment
functional $\L$ is called 
semi-classical if  there exist two (minimal) fixed polynomials $A(x)$ and
$B(x)$ with the properties that 
\be
\L\le(-B(x)p'(x) + A(x)p(x)\ri) = 0, \ \forall p(x)\in \C[x]\ .
\ee
The integral representation was obtained in \cite{ismail,marce,marce1}:
we can quickly reprove here their result (without details) in a different way which
was not used there and which is in the line of approach of this paper.
Consequence of the definition is that  the (possibly
formal) generating power series
\be
F(z):= \sum_{k=0}^{\infty} \mu_k \frac {z^k}{k!} \le(''='' \L({\rm
e}^{xz}) \ri) \ ,\ \
\mu_k:=\L(x^k)\ ,
\ee
satisfies the $n$-th order ODE
\be
\le[zB\le(\frac d{dz}\ri) - A\le(\frac d{dz}\ri)\ri]F(z)= 0\ .\label{genfu}
\ee
The order $n$ is the highest of the degrees of $A(x), B(x)$ and it is
referred to --in this context-- as  the {\em class}.
A distinction occurs according to the cases $\deg(A)< \deg
B$ (Case A in \cite{marce}) or $\deg(A)\geq \deg(B)$ (Case B).
By looking at the recursion relation satisfied by the moments $\mu_k$
one realizes that there are precisely $n$ linearly independent
solutions if in Case B or $n-1$ in Case A\footnote{In Case A and if
$A(x)\not\equiv 0$ there is
a linear constraint on the initial conditions for the recurrence
relation, which decreases the dimension of solution space by one. If
$A(x)\equiv 0$ then the solutions of the functional equation can be
found easily.} 
 and hence the functionals are in one--to--one correspondence
with the solutions of Eq. (\ref{genfu}) which are analytic at $z=0$.\\
It is precisely the result of \cite{miller} that the fundamental
system of solutions of Eq. (\ref{genfu}) are expressible as Laplace
integral transform  of the weight density 
\be
W(x):= \exp\le(\int {\rm d x} \frac {A(x)+B'(x)}{B(x)} \ri)\ ,
\ee
(which may have also branch-points) over $n$ distinct suitably chosen
contours $\Gamma_j$;
\be
F_j(z):= \int_{\Gamma_j} {\rm d}x W(x){\rm e}^{xz}\ .
\ee
In Case A one should actually reject one solution among them,
i.e. the one with a singularity at the origin, or better consider only
the linear combinations which are analytic at $z=0$.\par
In the present paper the bi-moment functionals we consider will rather 
correspond to generating functions in two variables satisfying an
over-determined (but compatible) system of PDEs, and the fundamental
solutions will be representable as suitably chosen double Laplace integrals.
The paper is organized as follows:\\
in Section \ref{definitions} we introduce the basic objects and
definitions, recalling how to explicitly construct the BOPs from the
matrix of bimoments. We also prove that the BOPs uniquely determine
the bi-moment functional: this is the analog in this setting of
Favard's Theorem which allows to reconstruct a moment functional from
any sequence of polynomials which satisfy a three--term recurrence
relation.\\
In Section \ref{bifunctionals} we introduce the definition of
semiclassical functionals and then prove that (under certain general
assumpions) they are representable as integrals of suitable 2-forms
over Cartesian products of complex paths.
The starting point is the fact already mentioned  that the generating
function of bi-moments now
depends on two variables $z,w$ and satisfies an over-determined
system of PDEs. We will prove the compatibility of this system (in the
class of cases specified in the text) and then we will solve it. The
solutions  
that we obtain (in the cases we consider) are entire functions of both
variables $z,w$ so that one could derive bounds on the growth of the
bi-moments (the coefficients of the Taylor series centered at $z=0=w$).\\
It should also be remarked that all semiclassical linear moment
functionals can be recovered as a special case of bilinear ones (see
Remark \ref{reduct}): this correspond to the fact that one-matrix
models can be recovered from two-matrix models in which one of the
measures is Gaussian.
\section{Definitions and first properties}
\label{definitions}
By bi-moment functional 
we mean a functional $\L$ on the tensor
product of two copies of the space of polynomials
\be
\L : \C[x]\otimes\C[y] \to \C\ .
\ee
Although the two polynomial spaces are just copies of the same space,
 we use two different indeterminates $x$ and $y$ in order to distinguish them.\\
Such a functional is uniquely determined by its
bi-moments
\be
\mu_{ij} := \L(x^i|y^j).
\ee
It makes sense  to look for bi-orthogonal polynomials.
We recall their definition and some standard facts \cite{mclaughlin,RM2}
\bd
Two sequences of polynomials $\{\pi_n(x)\}_{n\in \N}$ and
$\{\sigma_n(y)\}_{n\in \N}$ of exact degree $n$ are said to be
biorthogonal with respect to the bi-moment functional $\mathcal L$ if 
\be
\L(\pi_n|\sigma_m) = \delta_{nm}\ .
\ee
\ed
If such two sequences exist then we denote by $\{p_n(x)\}_{n\in \N}$ and
$\{s_n(y)\}_{n\in \N}$  the corresponding sequences of monic
polynomials, which then satisfy
\be
\L(p_n|s_m) = h_n\delta_{nm}\ ,\qquad h_n\neq 0 \ ,\forall n\in \N.
\ee
It is an adaptation of the classical result for orthogonal
polynomials to write a formula for the monic sequences
\bp
The biorthogonal polynomials exist if and only if
\be
\Delta_n\neq 0,\ n\in\N,\qquad 
 \Delta_n:= \det\pmatrix{
\!\!\mu_{0,0}\!\! & \!\!\mu_{0,1}\!\! & \!\!\cdots\!\! &\!\! \mu_{0,n\!-\!1}\!\! \cr
\!\!\mu_{1,0}\!\! &\!\! \mu_{1,1}\!\! &\!\! \cdots\!\! &\!\! \mu_{1,n\!-\!1}\!\! \cr
\!\!\vdots \!\! &\!\! \cdots\!\!  & \!\!\cdots\!\! &\!\! \vdots \!\!  \cr
\!\!\mu_{n\!-\!1,0}\!\! & \!\!\mu_{n\!-\!1,1}\!\! & \!\!\cdots \!\!&\!
\!\mu_{n\!-\!1,n\!-\!1}\!\!}\ , 
\ee
Under this hypothesis 
the monic sequences $\{p_n\}_{n\in \N}$ and $\{s_n\}_{n\in \N}$  
are given by the formulas
\bea
&& p_n(x):=\frac 1{\Delta_n} \det\pmatrix{
\!\mu_{0,\!0} \! &\!\! \cdots\!\! &\!\! \mu_{0,n-1}\!\! &\!\! 1\cr\vspace{-3pt}
\!\mu_{1,0}\!\!  & \!\!\cdots\!\! &\! \!\mu_{1,n-1} \!\!&\!\! x\cr\vspace{-3pt}
\!\vdots\!\!  & \!\!\cdots\!\! &\!\!\cdots\!\! &\!\! \vdots \! \! \cr\vspace{-3pt}
\!\mu_{n,0}\!\! & \!\!\cdots\!\! &\!\! \mu_{n,n-1}\!\! & \!\!x^{n}\!}\ ;\\
&& s_n(y):=\frac 1{\Delta_n} \det\pmatrix{
 \!\mu_{0,0} \! \! & \!\!\cdots\!\! & \!\!\mu_{0,n-1}\!\! &\!\! \mu_{0,n}\!\!\cr
\!\mu_{1,0}\!\! &\!\! \cdots\!\! & \!\!\mu_{1,n-1}\!\! &\!\! \mu_{1,n}\!\!\cr
\!\vdots \!\!   & \!\!\cdots\!\! &\!\!\cdots\!\! &\!\! \vdots\! \!  \cr
\!1\!\! & \!\!\cdots \!\!&\!\!y^{n-1}\!\! &\!\! y^{n}\!} \ .\label{biortho}
\eea
\ep
The proof of this simple proposition is essentially the same as for
the orthogonal polynomials and it is left to the reader (see
\cite{RM2,mclaughlin}).\\
With formula (\ref{biortho}) we can also compute 
\be
\L(p_n|s_m) = \frac{\Delta_{n+1}}{{\Delta_n}}\delta_{nm}\ .
\ee
The relation with the normalized polynomials is 
\be
\pi_n(x) = c_n p_n(x)\ ;\ \ \sigma_n(y):= \tilde c_n s_n(y)\ ,
\ee
where the complex constants $c_n$ and $\tilde c_n$ are such that
$c_n\tilde c_n =\frac { \Delta_{n+1}}{\Delta_n}$. \\
If biorthogonal polynomials exist they in general do not satisfy a
three terms recurrence relation as for the ordinary orthogonal
polynomials: they rather satisfy  recurrence relations which
generally are  not of finite bands
\bea
x\pi_n(x) = \gamma_n\pi_{n+1}(x) + \sum_{j=0}^n a_j(n)\pi_{n-j}(x)
\label{recrels1}\\
y\sigma_n(y) =\tilde \gamma_n\sigma_{n+1}(y) \sum_{j=0}^n
b_j(n)\sigma_{n-j}(y)\ .
\label{recrels2}
\eea
In the case of orthogonal polynomials the three terms recurrence
relation is sufficient for reconstructing the moment functional
(Favard's Theorem \cite{chihara}). A natural question is whether the
recurrence relations (\ref{recrels1}, \ref{recrels2}) are also sufficient for the
existence of a moment bifunctional for which the two sequences are
bi-orthogonal polynomials. Note that the specification of
the numbers $\gamma_n,\alpha_i(n)$, $i\leq n$ and  $\tilde \gamma_n,\beta_i(n)$,
$i\leq n$ determines uniquely  the two sequences of polynomials (with the
understanding that $\pi_{-n}\equiv 0\equiv \sigma_{-n}$)  in
Eqs. (\ref{recrels1},\ref{recrels2}) provided that
 $\gamma_n\neq 0\neq \tilde \gamma_n$,
$\forall n\in\N$. 
The following theorem answers 
 positively to the existence of the moment bifunctional
\bt[Favard-like Theorem for biorthogonal polynomials]
If the constants $\gamma_n, \tilde \gamma_n$ do not vanish for all $n\in \N$
then there exists a unique moment bifunctional $\mathcal L$ for which
the two sequences of polynomials $\pi_n,\sigma_n$ as in
Eq. (\ref{recrels1}, \ref{recrels2}) are biorthogonal.
\et
{\bf Proof}.
As for the ordinary Favard's theorem we proceed to the
construction of the bi-moments $\mu_{ij} = \mathcal L(x^i|y^j)$ by
induction. We introduce the associated  monic polynomials by defining 
\bea
&& p_n(x) := \frac 1{\pi_0} \pi_n(x)\prod_{k=0}^{n-1} \gamma_k \ ,\qquad
p_0(x)\equiv 1,\\
&& s_n(y) :=\frac 1{\sigma_0}\sigma_n(y) \prod_{k=0}^{n-1} \tilde
\gamma_k\ ,\qquad s_0(y)\equiv 1\ .
\eea
 The corresponding recurrence relations have the same form as in
Eq. (\ref{recrels1}, \ref{recrels2}) except that now the constants $\gamma_n, \tilde \gamma_n
$ are replaced by $1$.\\
The first moment $\mu_{00}$ is fixed by the requirement
\be
1 =  \mathcal L(\pi_0|\sigma_0)=\mu_{00} \pi_0\sigma_0\ ,
\ee
since the polynomials $\pi_0, \ \sigma_0$ are just nonzero constants.\\
Suppose now that the moments $\mu_{ij}$ have already been defined for
$i,j<N$. We need then to define the moments $\mu_{Nj}$ for $j=0,\dots
N-1$, and $\mu_{iN}$ for $i=0,\dots,N-1$ and $\mu_{NN}$.
By imposing the orthogonality
\be
0=\L(p_N|s_0) = \mu_{N0} + \dots\ ,
\ee
we define $\mu_{N0}$, where the dots represent an expression which contains
only moments already defined (i.e. $\mu_{i0}$, $i<N$).
We define by induction on $j$ the moments $\mu_{Nj}$, the first having
been defined above. We have, for $j<N-1$
\be
0=\L(p_N|s_{j+1}) = \mu_{N,j+1}+\dots\ ,
\ee
where again the dots is an expression involving only previously
defined moments. This defines $\mu_{N,j+1}$. We can repeat the arguments
for the moments $\mu_{iN}$, $i<N$ by reversing the role of the $p_i$'s and
$s_j$'s.\\
Finally the moment $\mu_{NN}$ is defined by 
\be
\det\pmatrix{\mu_{00} & \cdots & \mu_{0N}\cr
\vdots & &\vdots\cr
\mu_{N0}&\cdots& \mu_{NN}} 
=\frac 1 {\pi_0\sigma_0}\prod_{k=0}^{N-1} \gamma_k\tilde\gamma_k\ , 
\ee
where the only unknown is precisely $\mu_{NN}$ and its coefficient in
the LHS does not vanish since the corresponding minor is just 
\be
\det\pmatrix{\mu_{00} & \cdots & \mu_{0N-1}\cr
\vdots & &\vdots\cr
\mu_{N-10}&\cdots& \mu_{N-1N-1}} 
=\frac 1 {\pi_0\sigma_0} \prod_{k=0}^{N-2} \gamma_k\tilde\gamma_k\neq 0\ . 
\ee
This completes the definition of the moment bifunctional
$\L$. Q.E.D. \par\vskip 3pt
We now turn our attention to some specific class of bilinear
functionals $\L$.
We do not require for the analysis to come that the biorthogonal
polynomials exist, although for applications to multimatrix models
this is essential. In those applications the determinants $\Delta_n$ are
proportional to  the partition functions for  the corresponding multi-matrix
integrals (up to a multiplicative factor of $n!$) and are also
interpretable as tau functions of KP and $2$-Toda hierarchies
\cite{AvM2,UT} 
\section{Bilinear semiclassical functionals} 
\label{bifunctionals}
The notion of semiclassical for ordinary moment functionals and the
applications to random matrices suggest
the following 
\bd
We say that a bilinear functional $\L:\C[x]\otimes_\C\C[y]\to \C$ 
 is {\em semiclassical} if there exist four polynomials
$A_1(x),B_1(x)$ and $A_2(y),B_2(y)$ of degrees
$a_1+1,b_1+1,a_2+1,b_2+1$
respectively, such that the following distributional equations are fulfilled
\bea
\le\{\begin{array}{l}
\ds{( D_x\circ B_1(x)   + A_1(x))\otimes \1 \L = B_1(x)\otimes y \L }\\[3pt]
\ds{\1\otimes(D_y\circ B_2(y) + A_2(y))\,\L =x\otimes B_2(y)\,\L}\ .
\end{array}\ri.
\eea
Explicitly these equations mean that, for any
polynomials $p(x),s(y)$
\bea
&& \L\Big(-B_1(x)p'(x)+ A_1(x) p(x)\Big|s(y)\Big) =
\L\Big(B_1(x)p(x)\Big|ys(y)\Big) \label{semi1}
\\
&&\L\Big(p(x)\Big|-B_2(y) s'(y) +A_2(y) s(y)\Big) =
\L\Big(xp(x)\Big|B_2(y)s(y)\Big)\label{semi2} 
\eea
\ed
\br
\label{reduct}
We mentioned that any semi-classical moment functional is --in a
certain sense-- a special case of bilinear semi-classical functional. We want to
clarify this relation here.
Let us consider a semiclassical bifunctional in which $A_2(y) = ay$
and $B_2(y)=1$. The defining relations become
\be
 \L(-B_1p'+ A_1p|s)= \L(B_1p|ys)\ ,\ \ \L(p|-s' + ays) = \L(xp|s).\label{330}
\ee 
In particular for $s(y)= 1$ the second in Eq. (\ref{330}) reads 
\be
\L(p|y) = \frac 1 a \L(xp|1)\ .
\ee
The claim that the reader can check directly is that the moment
functional $\L_r(\cdot):= \L(\cdot|1)$ is a semiclassical
functional in the sense explained in the introduction with $A(x) = A_1(x) -\frac
x{a} B_1 (x)$ and $B(x) = B_1(x)$. It will be clear later on that this
``reduction'' corresponds to a partial integration of a Gaussian weight.
\er
In analogy with the orthogonal polynomials case we also define the
class
\bd
For a semi-classical bi-functional $\L$ we define its {\em bi-class}
as the pair of integers
\be
(s_1,s_2) = (\max(a_1,b_1)+1,\max(a_2,b_2)+1)\ .
\ee
\ed
Note that from the definition some recurrence relations follow for
the moments $\mu_{ij}$. In order to spell them out we
introduce the following notations for the coefficients of the
polynomials $A_i,B_i$
\bea
&& A_1(x)  = \sum_{j=0}^{a_1+1} \alpha_1(j)x^j\ ;\ \ B_1(x):=
\sum_{j=0}^{b_1+1} \beta_1(j)x^j\\
&&  A_2(y)  = \sum_{j=0}^{a_2+1} \alpha_2(j)y^j\ ;\ \ B_2(y):=
\sum_{j=0}^{b_2+1} \beta_2(j)y^j\ .
\eea
Then the aforementioned recurrence relations are given by 
\bp
The moments $\mu_{ij}$ of the classical bi-functional $\L$ are subject
to the relations
\bea
&& \sum_{j=0}^{b_1+1} \beta_1(j)\mu_{n+j,m+1} =- n\sum_{j=0}^{b_1+1}
\beta_1(j) \mu_{n-1+j,m}  + \sum_{j=0}^{a_1+1}\alpha_1(j)
\mu_{n+j,m}\label{momrec1} \\
&& \sum_{j=0}^{b_2+1} \beta_2(j) \mu_{n+1,m+j} =
-m\sum_{j=0}^{b_2+1}\beta_2(j) \mu_{n,m-1+j} + \sum_{j=0}^{a_2+1}
\alpha_2(j) \mu_{n,m+j}\ .\label{momrec2}
\eea
\ep
{\bf Proof.}\\
From the definition of semi-classicity by setting
$p(x)= x^n$ and $s(y)=y^m$ in the two relations (\ref{semi1},
\ref{semi2}). Q.E.D.\par\vskip 3pt
The two recurrence relations give an overdetermined system for the
moments: it is not guaranteed a priori that solutions exist and if
they do, how many.
There are now four
different cases, according to $\deg(B_i)\ds{\m{=}^{<}_{>}}
\deg(A_i)$; 
we  address in the present paper the case $\deg(A_i)>\deg(B_i)$,
$i=1,2$ (most relevant in the
applications to random matrix models) which is the analog of Case B
in \cite{marce} and we could call ``Case BB''. 
The other cases have less interesting applications in matrix models
because they correspond to potentials (in a sense which will be clear
below) which are bounded at infinity. They are certainly
interesting from the point of view
of Eqs. (\ref{momrec1}, \ref{momrec2});
 for example it is a simple exercise to check that if
$\deg(B_1) = \deg(B_2)=1 $ and $\deg(A_1)=\deg(A_2)=0$ then in general
no nontrivial solutions exist for Eqs  (\ref{momrec1}, \ref{momrec2}).\par
For the rest of this paper we will make the
following\\
{\bf Assumptions ($\mathcal A$)}
\be
\deg(B_i)+1 \leq \deg(A_i)\ ,\qquad i=1,2.
\ee
Moreover in the case
$\deg(B_1)+1=\deg(A_1)$ and $\deg(B_2)+1=\deg(A_2)$ we impose 
\be
\det\pmatrix{ \alpha_1(a_1+1)& \beta_1(b_1+1)\cr \beta_2(b_2+1)
& \alpha_2(a_2+1)}\neq 0 \ \ \hbox{when } a_1=b_1+1,\ a_2=b_2+1\ .\label{determin}
\ee
Under this assumption we can prove 
\bp
\label{dimension}
The solutions to Eqs. (\ref{momrec1}, \ref{momrec2}) form a vector space of
dimension $M:= s_1\cdot s_2 = (a_1+1)\cdot(a_2+1)$.
\ep
{\bf Proof.} The fact that the space of solutions is a vector space is
obvious from the linearity of the defining equations.
 We need to prove the assertion regarding
the dimension.\\
We define the (possibly formal) generating function of moments
\be
F(z,w):= \sum_{j,k=0}^\infty \frac {z^j w^k}{j!k!} \mu_{jk} =
\L\Big({\rm e}^{xz}\big|{\rm e}^{yw}\Big)\ .
\ee
From the recursion relation for the moments or
(equivalently) from the definition of semi-classicity, it follows that
such function satisfies the system of  PDEs
\bea
\le\{\begin{array}{l}
\displaystyle{\Big[(\pa_z+w) B_2(\pa_w) - A_2(\pa_w)\Big]F(z,w)=0}\\[10pt]
\displaystyle{\Big[(\pa_w+z) B_1(\pa_z\,) - A_1(\,\pa_z\,)\Big]F(z,w) = 0}
\end{array}\ri.\label{oversys}
\eea
Conversely, any solution of this system which is analytic at $z=0=w$ 
 provides a semi-classical bi-moment functional associated with the
data $A_i,B_i$. We now count the solutions of this system. It will be
clear later on that all the solutions are analytic at $z=0=w$ (in fact
entire) so that
any solution does define a moment functional.\par 
The system (\ref{oversys}) is a higher order overdetermined 
system of PDEs for the single function (or
formal power series) $F(z,w)$ and the compatibility is readily seen since 
\bea
&& \Big[(\pa_z+w) B_2(\pa_w)
-A_2(\pa_w),(\pa_w+z)B_1(\pa_z)  - A_1(\pa_z)
\Big] =\\
&&= \Big[(\pa_z+w) B_2(\pa_w),(\pa_w+z)B_1(\pa_z) \Big]  = \\
&&= \Big[(\pa_z+w),(\pa_w+z)\Big] B_2(\pa_w)B_1(\pa_z) = (1-1)
B_2(\pa_w)B_1(\pa_z) =0.
\eea
Now we express the system as a first order system of PDE's on the
suitable jet extension. Let us
introduce the notation 
\be
F_{\mu,\nu}(z,w):= {\pa_z}^\mu {\pa_w}^\nu F(z,w).
\ee
The proof now proceeds according to the three different cases: 
\begin{enumerate}
\item []Case BB1: $\deg(A_i)\geq \deg (B_i)+1$, $i=1,2$;
\item []Case BB2: $\deg(A_1)=\deg(B_1)+1$ but $\deg(A_2)>\deg(B_2)+1$ (or
vice-versa);
\item []Case BB3: $\deg(A_1)=\deg(B_1)+1$, $\deg(A_2)=\deg(B_2)+1$.
\end{enumerate}
For convenience 
we set the leading coefficients of the two polynomials $A_i$ to unity
as this does not affect the dimension of the solution space of the
system but make the formulas to come shorter to write.\\
In Case BB1  ($a_i\geq b_i+2$) we can write the two first
order systems for the  systems 
{\scriptsize
\bea
\le\{\begin{array}{lr}
\ds{ \pa_w F_{\mu,\nu} = F_{\mu,\nu+1}}  &\hspace{-5cm} \mu=0,\cdots,a_1,\
\nu=0,\cdots  a_2-1\\
\ds{\pa_w F_{\mu,a_2} = \sum_{k=0}^{b_2+1}\beta_2(k)(
wF_{\mu,k} + F_{\mu+1,k})
-\sum_{k=0}^{a_2} \a_2(k)F_{\mu,k}}  & \hspace{-5cm}  \mu=0..a_1-1\\
\ds{\pa_w F_{a_1,a_2} =  \sum_{k=0}^{b_2+1}\beta_2(k)\le[
wF_{a_1,k} + \le( 
\sum_{j=0}^{b_1+1} \b_1(j)\Big(zF_{j,k}+ F_{j,k+1} \Big)
-\sum_{j=0}^{a_1} \a_1(j)F_{j,k}
\ri) 
\ri]
-\sum_{k=0}^{a_2} \a_2(k)F_{a_1,k}}&
\end{array}\ri.
\eea
\bea
\le\{\begin{array}{lr}
\ds{ \pa_z F_{\mu,\nu} = F_{\mu+1,\nu}}  &\hspace{-5cm} \mu=0,\cdots,a_1-1,\
\nu=0,\cdots  a_2
\\
\ds{\pa_z F_{a_1,\nu} = \sum_{j=0}^{b_1+1}\beta_1(j)(
zF_{j,\nu}  +  F_{j,\nu+1})
-\sum_{j=0}^{a_1} \a_1(j)F_{j,\nu}}  & \hspace{-5cm}  \nu=0..a_2-1
\\
\ds{\pa_z F_{a_1,a_2} =  \sum_{j=0}^{b_1+1}\beta_1(j)\le[
zF_{j,a_2}  + \le( 
\sum_{k=0}^{b_2+1} \b_2(k)\Big(wF_{j,k}+ F_{j+1,k} \Big)
-\sum_{k=0}^{a_2} \a_2(k)F_{j,k}
\ri) 
\ri]
-\sum_{j=0}^{a_1} \a_1(j)F_{j,a_2}}
&
\end{array}\ri.\label{uno}
\eea}
Note that the two systems are consistent for the unknowns
$F_{\mu,\nu}$, $\mu=0,...,a_1$, $\nu=0,...,a_2$ if we have 
$b_i+2\leq a_i$, $i=1,2$.\\
In Case BB2 with $a_1=b_1+1$ the second system is not anymore
consistent because the RHS of the third equation in system (\ref{uno})
contains $F_{a_1+1,a_2}$. It must be replaced by 
{\scriptsize
\bea
\le\{\begin{array}{lr}
\ds{ \pa_z F_{\mu,\nu} = F_{\mu+1,\nu}}  &\hspace{-5cm} \mu=0,\cdots,a_1-1,\
\nu=0,\cdots  a_2
\\
\ds{\pa_z F_{a_1,\nu} = \sum_{j=0}^{a_1}\big(\beta_1(j)(
zF_{j,\nu}  +  F_{j,\nu+1})-\a_1(j)F_{j,\nu}\big)}  & \hspace{-5cm}  \nu=0..a_2-1
\\
\ds{\pa_z F_{a_1,a_2} = 
 \sum_{j=0}^{a_1}\beta_1(j)\le[
zF_{j,a_2}  + \le( 
\sum_{k=0}^{b_2+1} \b_2(k)wF_{j,k}
-\sum_{k=0}^{a_2} \a_2(k)F_{j,k}
\ri) 
\ri]
-\sum_{j=0}^{a_1} \a_1(j)F_{j,a_2}+}& \\
\ds{ +  \sum_{j=0}^{a_1-1}
\sum_{k=0}^{b_2+1} \beta_2(k)\beta_1(j) F_{j+1,k} + \beta_1(a_1)
\sum_{k=0}^{b_2+1} \beta_2(k)\le(\sum_{j=0}^{a_1}\bigg(\beta_1(j)(
zF_{j,k}  +  F_{j,k+1}) 
- \a_1(j)F_{j,k}\bigg)\ri) }&
\end{array}\ri.
\eea}
Finally in the Case BB3 ($a_1=b_1+1$ and $a_2=b_2+1$) we have the two
systems 
{\scriptsize
\bea
\le\{\begin{array}{lr}
\ds{ \pa_z F_{\mu,\nu} = F_{\mu+1,\nu}}  &\hspace{-5cm} \mu=0,\cdots,a_1-1,\
\nu=0,\cdots  a_2
\\
\ds{\pa_z F_{a_1,\nu} = \sum_{j=0}^{a_1}\big(\beta_1(j)(
zF_{j,\nu}  +  F_{j,\nu+1})-\a_1(j)F_{j,\nu}\big)}  & \hspace{-5cm}  \nu=0..a_2-1
\\
\ds{(1-\beta_1(a_1)\beta_2(a_2)) \pa_z F_{a_1,a_2} = 
\sum_{j=0}^{a_1} \beta_1(j)\le[z F_{j,a_2} + \sum_{k=0}^{a_2}
\big(w\beta_2(k)-\alpha_2(k) \big) F_{j,k}\ri] +}\\
\ds{-
\sum_{j=0}^{a_1}\alpha_1(j) F_{j,a_2} +
\sum_{j=0}^{a_1-1}\sum_{k=0}^{a_2}\beta_1(j)\beta_2(k)\pa_z F_{j,k}
 }& 
\end{array}\ri.
\eea}
and a similar system for the $\pa_w$ derivative. Note that in the
third equation the derivatives $\pa_z F_{j,k}$ are defined by the first
and second equation.\\
Since now 
 $(1-\beta_1(a_1)\beta_2(a_2))\neq 0$ as per the {\bf Assumption} (which is
$(\alpha_1(a_1+1)\alpha_2(a_2+2) -\beta_1(a_1)\beta_2(a_2))\neq 0$ if
we do not assume that the polynomials $A_1,A_2$ are monic) then the
system is still well defined; on the other hand, if
$(1-\beta_1(a_1)\beta_2(a_2))=0$ then the last equation becomes a {\em
constraint}\footnote{We are not going to examine this case in this
paper because it is more natural to study in the context of
semiclassical functionals of type AB or AA, i.e. when $\deg(A_i)\leq
\deg (B_i)$}.\par
It is a lengthy but straightforward check that the two systems are
indeed compatible in each of the three cases. Since the size of the system is
$M=(a_1+1)\cdot(a_2+1)=s_1s_2$ then there are precisely $M$ linearly
independend solutions. Q.E.D.\par\vskip 3pt
\br
In principle we would not have to check the compatibility because we
will construct later $M=s_1s_2$ solutions to the system, which
therefore will be proven to be  compatible {\em a posteriori}:
 the point of Prop. \ref{dimension} is
principally that the dimension of the solution space certainly does
not exceed $M$ because that is the dimension of a closed system in
the jet space.
\er
The Proposition implies that the recurrence relations (\ref{momrec1},
\ref{momrec2}) determine uniquely the functional $\L$ in terms of the
moment $\mu_{ij}$ with $i=0,\dots,a_1,j=0,\dots,a_2$. 
We need to produce  $M=s_1s_2$
linearly independent semiclassical functionals associated to the same
data $(A_1,B_1,A_2,B_2)$ by means of integral representations.\\
Equivalently we can produce integral representation for the $M$
linearly independent solutions of the overdetermined system of PDE's
(\ref{oversys}).
It is precisely in this form that we will solve the problem, showing
contextually that the generating functions are indeed entire functions
of $w,z$.
The starting point is to assume that such an integral representation
exists:  so suppose that 
\be
F(z,w) = \int_{\Gamma_{(x)}}\int_{\Gamma_{(y)}}{\rm
d}x\wedge {\rm d}y W(x,y) {\rm e}^{xz+yw}\ ,   
\ee
is a double Laplace integral representation for a solution of
(\ref{oversys})\footnote{In principle one could integrate the two-form
$W(x,y){\rm e}^{xz+yw}{\rm d}x\wedge {\rm d}y$ over any $2$-cycle, but
here we do not need such generality}.\\
Plugging such representation in the two equations  in (\ref{oversys})
and assuming that the contours are so chosen as to allow integration
by parts without boundary terms, we obtain two first order equations
for the bi-weight $W(x,y)$ 
\bea
\Big(B_1(x)\pa_x + A_1(x) + B_1'(x)\Big) W(x,y) = y\,B_1(x) W(x,y)\label{qweq1}\\
\Big(B_2(y)\pa_y + A_2(y) + B_2'(y)\Big) W(x,y) = x\,B_2(y) W(x,y)\ .\label{qweq2}
\eea
We make the {\bf Assumption ($\mathcal B$)}  
 that each pair $(A_i,B_i)$ are relatively prime or at
most share a factor $(x-c)$ (or $(y-s)$).
 The reason is similar to the case of ordinary semiclassical
functionals.  We will return on this genericity assumption later on.\par
The two differential equations (\ref{qweq1},\ref{qweq2}) 
form an overdetermined system for the biweight $W(x,y)$ which is
compatible and can be solved to give the
 only solution (up to a multiplicative nonzero constant)
\bea
&&W(x,y) =  W_1(x) W_2(y) {\rm e}^{xy} = \exp\le(-V_1(x) - V_2(x) +
xy\ri)\ ,\\[10pt] 
&& \frac {W_1'(x)}{W_1(x)} =  \frac {A_1(x)+B_1'(x)}{B_1(x)},\qquad
 \frac {W_2'(y)}{W_2(y)} =  \frac {A_2(y)+B_2'(y)}{B_2(y)}\ ,\label{logder}\\
&&V_1(x):=\int {\rm d}x \frac {A_1(x)+B_1'(x)}{B_1(x)} \\
&&V_2(y):= \int{\rm d}y \frac {A_2(y)+B_2'(y)}{B_2(y)}\ .
\eea
We call the two functions $V_1(x)$, $V_2(y)$ the {\em potentials}
(borrowing the name from the statistical mechanic and random matrix context).\\
Note that if there are nonzero residues at the poles of
$\frac{A_i+B_i'}{B_i}$ then the corresponding potential have
logarithmic singularities or poles.
The general form of the biweight is  
\bea
 W_1(x):= \prod_{j=1}^{p_1}(x-X_j)^{\lambda_j} \exp\le[V_1^+(x)  +
\frac {M_1(x)}{\prod_{j=1}^{p_1}(x-X_j)^{g_j}}\ri]\ ,\label{weightx} \\
 \deg(M_1)\leq
\sum_{j=1}^{p_1} g_j\ ,\ M_1(X_j)\neq 0 \nonumber\\
 W_2(y):= \prod_{k=1}^{p_2}(y-Y_j)^{\rho_k} \exp\le[V_2^+(y)  +
\frac {M_2(y)}{\prod_{k=1}^{p_2}(y-Y_k)^{h_k}}\ri]\ ,\label{biweight} 
\\ \deg(M_2)\leq
\sum_{k=1}^{p_2} h_k\ ,\ M_2(Y_k)\neq 0 \nonumber\ .
\eea
In this formulas and in the rest of the paper $X_j$ denote the zeroes of
$B_1(x)$, $g_j+1$ the corresponding multiplicities and $-\lambda_j$
are the residues at $X_j$ of the differential ${\rm d}V_1(x)$;
similarly, $Y_k$ denote the zeroes of $B_2(y)$, $h_k+1$ the corresponding
multiplicities and $-\rho_k$ the residues at $Y_k$ of the differential
${\rm d}V_2(y)$.\par
The bi-class of the corresponding semiclassical bifunctional is then
the total degree of the divisor of poles of the derivatives of the two
potentials on the Riemann
spheres whose affine coordinates are $x$ and $y$
\be
s_1 = d_1 + \sum_{j=1}^{p_1} (g_j+1)\ ,\ \ s_2 = d_2 + \sum_{j=1}^{p_2}
(h_j+1) \ .
\ee
We will also use the notations  $X_0=\infty\in {\mathbb
P}^1_x$, $Y_0=\infty\in {\mathbb P}^1_y$.
\subsection{The functionals}
We will define two sets of paths in the two punctured Riemann spheres
${\mathbb P}^1_x$ and ${\mathbb P}^1_y$. We focus on the
first sphere, the paths in the second being defined in analogous way.\\
More precisely we define $s_1$ ``homologically'' independent
paths in  ${\mathbb P}^1_x\setminus C_x$ and $s_2$
paths in ${\mathbb P}^2_y\setminus C_y$ where $C_x$ and $C_y$ are
suitable union of cuts and points: for example the set $C_x$ is the
union of all poles and essential singularities of $W_1(x)$ and cuts
extending from the branchpoints to infinity.\\
The reference to the homology is not in the ordinary sense: here 
 we are considering in fact the {\em relative} homology of the
cut-punctured sphere with prescribed sectors around the punctures.\\
We first define some sectors $S_k^{(j)}$, $j=1,\dots p_1$, $k=0,\dots g_j-1$.
around the points $X_j$ for which $g_j>0$
(the multiple zeroes of $B_1(x)$) in such a way that 
\be
\Re\le(V_1(x)\ri) \mathop{\longrightarrow}_{
\shortstack{{\scriptsize $x\to X_j,$}\\
{\scriptsize  $x\in S
^{(j)}_k$}}}
+\infty\ .
\ee
The number of sectors for each pole is the degree of that pole in the
exponential part of $W_1(x)$, that
is $d_1+1$ for the pole at infinity and $g_j$ for the $j$-th pole.
Explicitly 
\bea
&& S^{(0)}_k := \le\{ x:\in \C;\ \ \frac {2k\pi -\frac \pi
2+\epsilon}{d_1+1}< \arg(x)+\frac{\arg(v_{d_1+1})}{d_1+1} < \frac {2k\pi +\frac \pi
2-\epsilon}{d_1+1} \ri\}\ ,\ \ k=0\dots d_1\ ;\label{sectinf}\\
&&  S^{(j)}_k := \le\{ x:\in \C;\ \ \frac {2k\pi -\frac \pi
2+\epsilon}{g_j}< \arg(x-X_j)+\frac{\arg(M_1(X_j))}{g_j} < \frac {2k\pi +\frac \pi
2-\epsilon}{g_j} \ri\}\ ,\\
&&\hspace{2cm} k=0,\dots,g_j-1,\ \ j=1,\dots,p_1\ .\nonumber
\eea
These sectors are defined precisely in such a way that approaching any
of the
essential singularities (i.e. an $X_j$ such that $g_j>0$) the function
$W_1(x)$ tends to zero faster than any power.\\ 
{\bf Definition of the contours}\\
The definition of the contours follows directly \cite{miller}, but
 we have to repeat it in both Riemann spheres. For the sake of
completeness we recall the way they are defined.
\begin{enumerate}
\item For any $X_j$ for which there is {\em no essential singularity}
(i.e. $g_j=0$), then we have two subcases
\begin{enumerate}
\item Corresponding to the $X_j$'s which are  branch points or a pole
 ($\lambda_j\in\C\setminus \N$), we  
take a loop starting at infinity in some fixed sector
$ S_{k_L}^{(0)}$  encircling the singularity and going back to infinity in the same
sector.
\item For the $X_j$'s which are regular points ($\lambda_j\in \N$) we
take a line joining
$X_j$ to infinity and approaching $\infty$ in the same sector $
S_{k_L}^{(0)}$ as before.
\end{enumerate}
\item For any $X_j$ for which there is an essential singularity
(i.e. for which $g_j>0$) we define $g_j$ contours starting from $X_j$
in the sector $S_0^{(j)}$ 
and returning to $X_j$ in the next (counterclockwise) sector. Finally
we join the singularity $X_j$ to $\infty$ by a path approaching $\infty$
within the sector $ S^{(0)}_{k_L}$ chosen at point 1(a).
\item For $X_0:=\infty$ we take $d_1$ contours starting at $X_0$ in tha
sector $S^{(0)}_k$ and returning at $X_0$ in the sector
$S^{(0)}_{k+1}$.\footnote{Note
that in our assumptions on the degrees of $A_i,B_i$ the degrees of the
essential singularity at infinity satisfy $d_1\geq 1\leq d_2$}. 
\end{enumerate}
For later convenience we also fix a sector  $\mathcal S_L$ of width
$\beta<\pi-\epsilon$ which contains the sector $S^{(0)}_{k_L}$ used
above. 
 The picture below gives an example of the typical situation, where
the light grey sector represents $\mathcal S_L$.
 We will make use also of the sector $\mathcal E$ which is
a sector within the dual sector\footnote{ We recall that 
for a given sector $\mathcal S$ centered around a ray ${\rm
arg}(z)=\alpha_0$ with width $A<\pi$, the {\bf dual sector} $\mathcal
S^\vee$ is the sector centered around the ray ${\rm arg}(z)=
\pi-\alpha_0$ and with width $\pi-A$.} of $\mathcal S_L$ (in dark shade of
grey in the picture): it is not difficult to realize that we can
always arrange contours in such a way that $\mathcal E$ is a small
sector below the real positive axis (if the leading coefficient of
$V_1^+$ is real and positive, otherwise the whole picture should be
rotated appropriately).\\
We shall also require that all contours do not intersect except possibly at
some $X_j$ and that each closed loop should either encircle only one
singularity or have one of the $X_j$ on its support.\par
The result of this procedure produces precisely $s_1$ contours. By
virtue of Cauchy's theorem the choice is largely arbitrary.\\
An important feature for what follows is that 
{\em when a contour $\Gamma_j$ is closed (on the sphere $\mathbb P^1_x$), then
$W_1(x)$ has a singularity and/or is unbounded in the region inside
$\Gamma_j$}. We will call this property the {\bf Property
($\mathbf \wp$)}.\\[5pt]
\centerline{
\parbox{8cm}{
\epsfxsize=8cm
\epsfysize=8cm
\epsffile{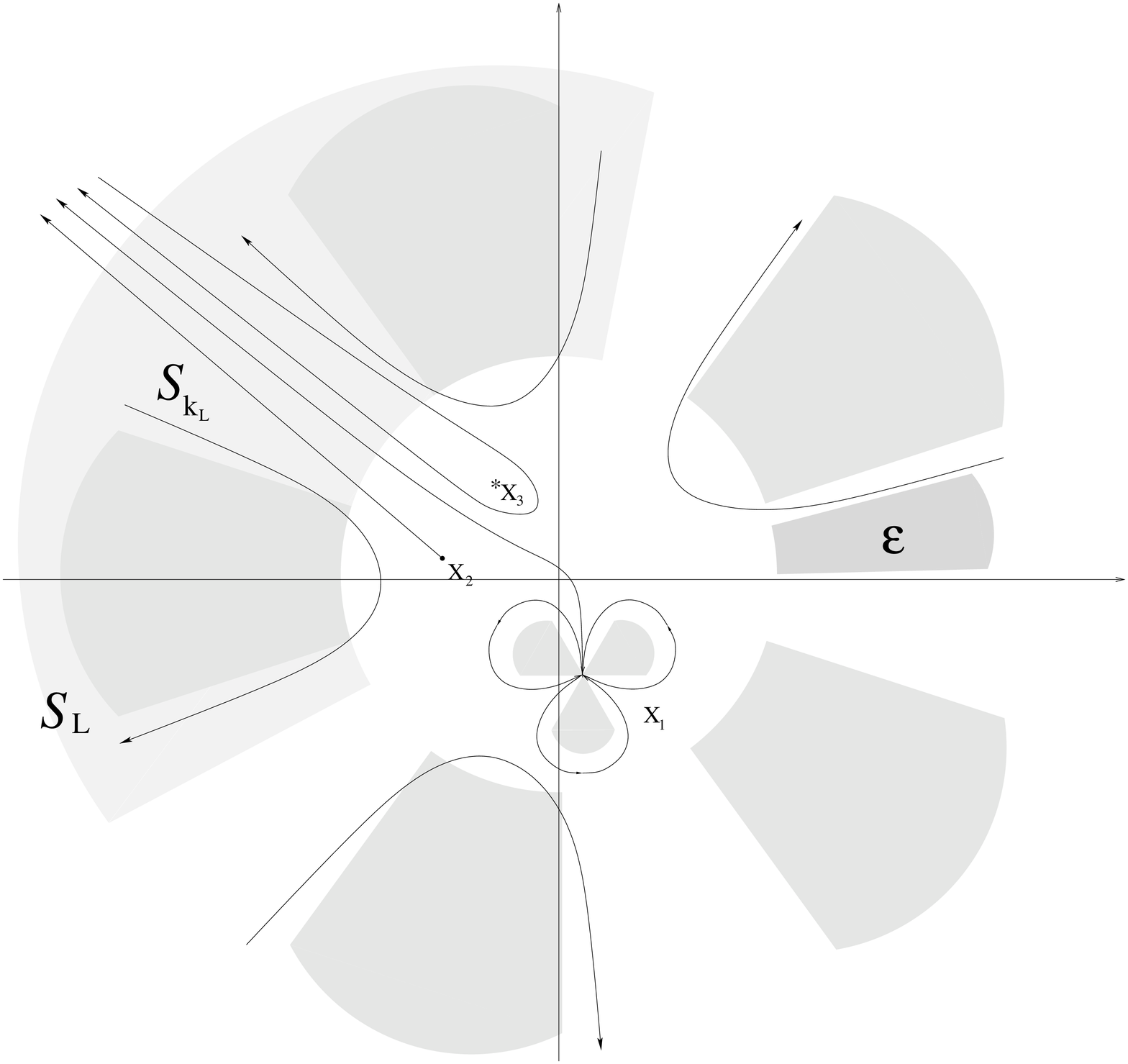}}}
\centerline{ \parbox{12cm}{Figure 1: The set of contours in the $x$ Riemann sphere
$\mathbb P^1_x$. Here we have three zeroes of $B(x)$, $X_1,X_2,X_3$,
and the singularity  at infinity
$X_0$ of order $d_1+1=5$.
 The zero $X_1$ has multiplicity $g_j+1=4$ and the corresponding
essential singularity behaves like $\exp{(x-X_1)^{-3}}$, the zero $ X_2$ is
a regular point for $W_1(x)$, namely $\lambda_2 \in \N$ and finally
the zero $X_3$ is either a branch point of $W_1$, in which case the
cut extends to infinity ``inside'' the contour (in the picture),  or a pole
($\lambda_3\not\in \N$).  
}}\vskip 3pt
We then
define the fundamental functionals by
\bea
\L_{ij} (x^n|y^m):= \int_{\Gamma^{(x)}_i\times \Gamma^{(y)}_j}
\!\! {\rm d}x\wedge {\rm d}y\,W_1(x)W_2(y){\rm
e}^{xy} x^n y^m
\ ,\ \\
 \ i=1,\dots,s_1,\ j=1,\dots s_2\ , n,m\in \N\nonumber . 
\eea
We point out that such contours are chosen so that the corresponding
functionals are defined on any monomials $x^j y^k$ and such that
integration by parts does not give any boundary contribution. 
Each such functional is a semi-classical functional
associated to the data $A_1,B_1,A_2,B_2$ and their number is precisely
the expected number $s_1s_2$ for the solutions of Eqs. (\ref{oversys}) for the
generating functions.
 The problem now is to show that they are linearly independent.
\br
A special care should be directed at the case $d_1=d_2=1$, i.e. when
$a_1=b_1+1$ and $a_2=b_2+1$. Indeed in this circumstance the two
polynomials $V_1^+(x)=\frac \delta 2 x^2 +..$ and $V_2^+(y)=\frac \sigma 2
y^2 + .. $ are just quadratic. The biweight
$W(x,y)$ has then the form 
\be
W(x,y) = \exp\le(-\frac \delta 2 x^2 -\frac \sigma 2 y^2 + xy + \dots
\ri)[\dots]\ .
\ee 
The condition on the determinant (\ref{determin}) is precisely the
nondegeneracy of the quadratic form $-\frac \delta 2 x^2 -\frac \sigma
2 y^2 + xy$. However, if $|\delta||\sigma|\leq 1$ then  the integrals as
we have defined are always divergent when two contours which stretch
to infinity are involved. This simply means that we cannot choose 
 the surface of integration in factorized form $\Gamma^{(x)}\times
\Gamma^{(y)}$ but need to resort to a surface which is not
factorized.\\
Alternatively  we can
analytically continue from the region of  $\delta,\sigma$ for which
the integrals are convergent. 
\er
Some important remarks are in order.
Consider the generating functions associated to these contours
\be
F_{ij}(z,w):= \int_{\Gamma^{(x)}_i\times \Gamma^{(y)}_j}
\!\! {\rm d}x\wedge {\rm d}y\,W_1(x)W_2(y){\rm
e}^{xy} 
 {\rm e}^{  xz + yw}\ .
\ee
They are entire functions of $z,w$ and
hence are indeed generating functions of the bi-moment functionals
$\mathcal L_{ij}(\cdot|\cdot)$. Indeed our assumptions on the
degrees guarantees that $V_i^+$ have degree at least $2$, which is
sufficient to guarantee analyticity w.r.t. $z,w$ in the whole complex plane.
\br
If the index $i$ corresponds to a bounded contour $\Gamma^{(x)}_i$ 
then $F_{ij}(z,w)$ is a
function of exponential type in $z$ (similarly for $w$ if
$\Gamma^{(y)}_j$ is bounded).
\er
\br If  the index $i$ corresponds to  one of the contours
$\Gamma^{(x)}_i$ defined at point 1(a)
or 1(b) above, then $F_{ij}(z,w)$ is of exponential type only for $z$
in an appropriate   sector which contains the sector $\mathcal E$
dual  to the sector $\mathcal S_L$.
\er
Before entering into the details of the proof of linear independence
let us return to the Assumption ($\mathbf{\mathcal B}$) about  the pairs
$(A_i,B_i)$.
 Suppose that -say- $A_1$ and $B_1$ have a common
factor $(x-c)^K$, $K\geq 1$ and that they have no other common
factor. That is let us suppose that
\bea
A_1(x) = (x-c)^l\tilde A_1(x)\ ,\ \ B_1(x) = (x-c)^r \tilde B_1(x)\ ,\
\\
 l>0<r,\ K:=\min(l,r)\ \nonumber ,
\eea
with $\tilde A_1(c)\neq 0\neq \tilde B_1(c)$.
 Then formula (\ref{logder}) would give 
\be
V'_1(x) = -\frac {W'_1(x)}{W_1(x)} = \frac {(x-c)^l\tilde A_1 + r(x-c)^{r-1}
\tilde B_1 + (x-c)^r \tilde B_1'}{(x-c)^r\tilde B_1} \ ,\label{symp}
\ee
so that the divisor of poles od ${\rm d}V_1(x)$ has degree {\em less}
than $s_1$.
Now we have two possible cases:\\
{\bf (i)} if $l\geq r-1$ then we can recast Eq. (\ref{symp}) in the
form 
\be
-\frac {W'_1(x)}{W_1(x)} = \frac {(x-c)^{l-r+1}\tilde A_1 +
(r-1)\tilde B_1 +\big((x-c)\tilde B_1\big)'}{(x-c)\tilde B_1} \ .\label{symp2}
\ee
which is equivalent to a problem in which the polynomials $A_1,B_1$
are substituted by $\underline{A_1}:= (x-c)^{l-r+1}\tilde A_1 +
(K-1)\tilde B_1$ and $\underline{B_1}:= (x-c)\tilde B_1$ respectively,
which now satisfy the assumption ({\bf F}). 
In particular the definition of the contours provides the correct
number of distinct contours for the new pair $
( \underline{A_1},\underline{B_1})$, that is
$s_1-r+1$ distinct contours (in the $x$ plane). We need to recover
$(K-1)s_2$ solutions if $l>r-1$ or $ls_2=Ks_2$ if $l=r-1$.\\
{\bf (ii)} If $l\leq r-2$ then  we can recast Eq. (\ref{symp}) in the
form 
\be
-\frac {W'_1(x)}{W_1(x)} = \frac {\tilde A_1 + l(x-c)^{r-1-l}
\tilde B_1 + \le((x-c)^{r-l} \tilde B_1\ri)'}{(x-c)^{r-l}\tilde B_1} \ .\label{symp3}
\ee
now equivalent to a problem in which the polynomials $A_1,B_1$
are substituted by $\underline{A_1}:= \tilde A_1 + K(x-c)^{r-l-1}\tilde B_1$ and 
$\underline{B_1}:= (x-c)^{r-l}\tilde B_1$ respectively,
which do not have the factor $(x-c)$ in common and hence satisfy the assumption ({\bf F}). 
The definition of the contours provides the correct
number of distinct contours for the new pair $
( \underline{A_1},\underline{B_1})$, and we need to recover $Ks_2$ solutions.\\
The next proposition shows how to recover the missing solutions.
\bp
If 
\be
A_1(x) = (x-c)^K\tilde A_1(x)\ ,\ \ B_1(x) = (x-c)^K \tilde B_1(x)\ ,\
\ K\geq 1\ ,
\ee
and $\tilde A_1(x)$, $\tilde B_1(x)$ do not vanish both at $c$ then
Eqs. (\ref{oversys}) have also the solutions
\be
F^{(j)}_k(z,w) = {\rm e}^{cz} \int_{\Gamma^{(y)}_k} \!\! {\rm d}y
(y+z)^j {\rm e}^{y(w+c)} W_2(y)\ ,\ \ j=0,...,K-1.\label{gendelta}
\ee
\ep
{\bf Proof.}\\
The fact that the functions (\ref{gendelta}) solve our system can be
checked directly.\\ 
Indeed the first eq. in (\ref{oversys}) is satisfied  because the
differential operator reads 
\be
(\pa_w + z)B_1(\pa_z) -A_1(\pa_z) = \le[(\pa_w + z)\tilde B_1(\pa_z)
-\tilde A_1(\pa_z)\ri](\pa_z-c)^K\ ,
\ee
and the proposed solutions are linear combination of functions 
of the form $z^r{\rm e}^{cz} f_r(w)$,\  $r<K$
which are all in the kernel of $(\pa_z-c)^K$.
The second equation in (\ref{oversys}) now reads
{\scriptsize
\beas
&& \le[ (\pa_z + w) B_2(\pa_w)-A_2(\pa_w)\ri] {\rm e}^{cz}
\int_{\Gamma^{(y)}_k} \!\! {\rm d}y
(y+z)^j {\rm e}^{y(w+c)} W_2(y) = \\
&&= c\,{\rm e}^{cz}
\int_{\Gamma^{(y)}_k} \!\! {\rm d}y\, B_2(y)
(y+z)^j {\rm e}^{y(w+c)} W_2(y)  + 
{\rm e}^{cz} \int_{\Gamma^{(y)}_k}\!\! {\rm d}y \Big(B_2(y)(\pa_z
+w) - A_2(y)\Big)  (y+z)^j {\rm e}^{y(w+c)} W_2(y) =\\
&&= {\rm e}^{cz}
\int_{\Gamma^{(y)}_k}\!\! {\rm d}y \Big ( B_2(y)(c+ \pa_z) -
A_2(y)\Big)
(y+z)^j {\rm e}^{y(w+c)} W_2(y) + {\rm e}^{cz}
\int_{\Gamma^{(y)}_k}\!\! {\rm d}y B_2(y)   W_2(y) (y+z)^j{\rm e}^{yc}
\pa_y({\rm e}^{yw})  = \\
&&=  {\rm e}^{cz}
\int_{\Gamma^{(y)}_k}\!\! {\rm d}y \Big(B_2(y) (\pa_z+c)  -
A_2(y)\Big)
(y+z)^j {\rm e}^{y(w+c)} W_2(y) + {\rm e}^{cz}
\int_{\Gamma^{(y)}_k}\!\! {\rm d}y B_2(y)   W_2(y) (y+z)^j
{\rm e}^{cy}\pa_y({\rm e}^{yw}) = \\
&&=  {\rm e}^{cz}
\int_{\Gamma^{(y)}_k}\!\! {\rm d}y \,W_2(y){\rm e}^{y(w+c)} 
\Big[B_2(y) \le[\pa_z-\pa_y\ri]\Big]
(y+z)^j   +\\
&&\hspace{1cm} +{\rm e}^{cz}
\int_{\Gamma^{(y)}_k}\!\! {\rm d}y  \Big(W_2'(y)B_2(y)-
(A_2(y) + B_1'(y))W_2(y)\Big)(y+z)^j{\rm e}^{y(w+c)} = 0.  
\eeas}
In Case {\bf(ii)} (or in Case {\bf(i)} but with $l=r-1$) 
these solutions are precisely the $Ks_2$ missing solutions.\\
In Case {\bf(i)} with $l\geq r$ only $l-1=K-1$ among the  solutions
(\ref{gendelta})  are
linearly independent from those defined in terms of the contour
integrals. To see this we write the
weight
\be
-\frac{W_1'(x)}{W_1(x)} = \frac r{x-c} + \frac {\tilde A_1 + \tilde
B_1'}{\tilde B_1}\ .
\ee
Since $\tilde B_1(c)\neq 0$ then $W_1(x)$ has a pole of order $r$ at
$x=c$ and can be written as 
\be
W_1(x)= (x-c)^{-r} w_1(x)\ ,
\ee
with $w_1(x)$ analytic at $x=c$ and $w_1(c)\neq 0$.
The contour which comes from infinity, encircles
$c$ and goes back to infinity can be retracted to a circle around the
pole, so that the corresponding solutions given by the integral
representation would be
\beas
\int_{\Gamma_y^{(k)}}\oint_{|x-c|=\epsilon}\!\!
 {\rm d}x\wedge {\rm
d}y \, (x-c)^{-r}w_1(x) {\rm e}^{x(z+y)+ wy} W_2(y) =\\
=
2i\pi(r-1)!\int_{\Gamma_y^{(k)}} \!\! {\rm d}y \pa_x^{r-1}\le.\le(w_1(x) {\rm
e}^{x(z+y)}\ri)\ri|_{x=c} W_2(y) \ .
\eeas
Such a solution is clearly an appropriate linear combination of the
$F^{(j)}_k$s $j=0,\dots r-1\leq K-1$
 with the nonzero coefficient $w_1(c)$ in front of
$F^{(r-1)}_k$. Q.E.D
\br
The function in Eq. (\ref{gendelta}) with $j=0$ corresponds to
a moment functional $\L = \delta_c\otimes \mathcal Y$, where $\mathcal
Y$ is any semi-classical moment functional associated to
$A_2(y),B_2(y)$ and $\delta_c$ is the delta functional  supported at
$x=c$ on the space of polynomials $\C[x]$.
The other solutions in Eq. (\ref{gendelta}) with
$j>0$ are also supported at $c$ but are not factorized and have the
form 
\be
\L = \sum_{k=0}^{j} \delta_{c}^{(k)}\otimes \mathcal Y_k\ ,
\ee
for suitable moment functionals $\mathcal Y_k$. 
\er
If there are other roots common to $A_i,B_i$ we can repeat the
procedure until we have a reduce 
problem which satisfies the Assumption ({\bf
$\mathcal B$}).\par
Therefore from this point on we will assume that the data
$(A_1,B_1,A_2,B_2)$ satisfy the Assumption ({\bf $\mathcal B$}).
\bt
\label{main}
The functionals $\L_{ij}$ or --equivalently-- the generating functions
\be
F_{ij}(z,w):= \int_{\Gamma^{(x)}_i\times \Gamma^{(y)}_j}
\!\! {\rm d}x\wedge {\rm d}y\,W_1(x)W_2(y){\rm
e}^{xy} 
 {\rm e}^{xz + yw}\ 
\ee
are linearly independent 
\et
The proof is an adaptation of \cite{miller} with a small improvement
(and a correction).
We prepare a few lemmas.
\bl [Theorem of Mergelyan (\cite{walsh}, p. 367)]
\label{margelyan}
If $E$ is a closed bounded set not separating the plane and if $F(z)$
is continuous on $E$ and analytic at the interior points of $E$, then
$F(z)$ can be uniformly approximated on $E$ by polynomials.
\el
The next Theorem is a rephrasing of the content of \cite{miller} for
the proof of which we refer ibidem.
\bt[Miller-Shapiro Theorem]
\label{millsha}
If $\Gamma$ is a closed simple Jordan curve and $F(z)$ is an analytic
function (possibly with singularities and/or multivalued) in the
points inside $\Gamma$ such that the equation  
\be
\oint_{\Gamma} F(z) p(z) {\rm d}z =0
\ee 
holds for any polynomial $p(z)\in (z-z_0)
\C[z]$ (for some fixed $z_0\in \Gamma$),
then $F(z)$ has no singularities inside $\Gamma$ and it is bounded in
the interior region of and on $\Gamma$.
\et
Suppose now by contradiction that there exist constants $C_{ij}$ not
all of which zero such that 
\be
\sum_{i=1}^{s_1}\sum_{j=1}^{s_2} C_{ij}\int_{\Gamma^{(x)}_i\times \Gamma^{(y)}_j}
\!\! {\rm d}x\wedge {\rm d}y\,W_1(x)W_2(y){\rm
e}^{xy} 
 {\rm e}^{  xz + yw}\equiv 0.\label{lindepcy}
\ee
{\bf Reduction of the problem}\\
We claim that if Eq. (\ref{lindepcy})
 holds then we also have
\be
0\equiv \sum_{i=1}^{s_1}\sum_{j=1}^{s_2} C_{ij}
\int_{\Gamma^{(x)}_i\times \Gamma^{(y)}_j}
\!\! {\rm d}x\wedge {\rm d}y\,W_1(x)W_2(y)
 {\rm e}^{  xz + yw} =  \sum_{i=1}^{s_1}\sum_{j=1}^{s_2} C_{ij}\, 
\Xi_i(z) \Psi_j(w)\ ,\label{lindepcyred}
\ee
where we have defined 
\bea
&&\Xi_i(z):= \int_{\Gamma^{(x)}_i}\!\! {\rm d}x \, W_1(x){\rm
e}^{xz} \\
&&\Psi_j(w):=  \int_{\Gamma^{(y)}_j}\!\! {\rm d}y \, W_2(y){\rm
e}^{yw}\ .
\eea
Indeed consider the auxiliary function of the new variable $\rho$ 
\be
A(\rho;z,w)
:= \sum_{i=1}^{s_1}\sum_{j=1}^{s_2} C_{ij}\int_{\Gamma^{(x)}_i\times \Gamma^{(y)}_j}
\!\! {\rm d}x\wedge {\rm d}y\,W_1(x)W_2(y){\rm
e}^{\rho xy + zx + wy}\ .
\ee
Here $z,w$ play the role of parameters.
This function is entire  in $\rho$ (because by our assumptions
$\deg(V_i^+)\geq 2$ and hence for all contours going to infinity the
integrand goes to zero at least as $\exp(-|x|^2-|y|^2)$), and by
applying $(\pa_z\pa_w)^K$ to
Eq. (\ref{lindepcy}) we have 
\be
0\equiv \le(\pa_z\pa_w\ri)^KA(1;z,w) = \le(\frac {d}{d\rho}\ri)^K A(\rho;z,w)
\bigg|_{\rho=1}\ , \
\ \forall K\in \N\ .
\ee
Therefore we also have $A(0;z,w)\equiv 0$, $\forall z,w\in \C$, which
is Eq. (\ref{lindepcyred}).\\
This shows that proving that the functions $F_{ij}$ are linearly
independent is
equivalent to proving that the two sets of  functions
$\{\Xi_i(z)\}_{i=1\dots s_1}$
and $\{ \Psi_j(w)\}_{j=1\dots s_2}$
are (separately) linearly independent.\\
Both the $\Xi_i$s and the $\Psi_j$s are now solutions of the decoupled
ODEs of the same type (i.e. with linear coefficients)
\bea
&&\le[zB_1\le(\frac d{dz}\ri)- A_1\le(\frac d{dz}\ri) \ri]\Xi_i(z) = 0\\
&&\le[wB_2\le(\frac d{dw}\ri)- A_2\le(\frac d{dw}\ri) \ri]\Psi_j(w) = 0\ .
\eea
Equivalently we may say that $\Xi_i$s and $\Psi_j$s are generating
functions for the moments of semiclassical functionals associated to
$(A_1,B_1)$ and $(A_2,B_2)$ respectively. Their linear independence 
 was proven in \cite{miller}. Unfortunately this latter paper has a
small flaw that makes one step of the  proof impossible when
$\deg(A_i)>\deg(B_i)+2$ (while it is correct if $\deg(A_i)\leq
\deg(B_i)+2$) \cite{shappriv}.\\
On the other side the linear independence of certain integral
representation for semi-classical moment functionals was obtained in
\cite{marce}; however their definitions for the contours forces them
to a procedure of regularization in certain cases which is elegantly
bypassed by the definition of the contours in \cite{miller}.
We prefer to fix the proof of \cite{miller} since then we will not
need any regularization.
\subsection{Linear independence of the $\Xi_i$s}
In this section we prove the linear independence of the 
 functions $\Xi_i$. This will also  prove 
the linear independence of the functions $\Psi_j$ since they are
precisely of the same form.
We assume that the polynomial $V_1^+(x)$ appearing in
Eq. (\ref{weightx})  has the form
\be
V_1^+(x) = \frac 1{d+1} x^{d+1} + \sum_{j=0}^{d} v_j x^j\ \
(d:=d_1\geq 1).\label{potent}
\ee
This does not affect the generality of the problem inasmuch as
it amounts to a rescaling of the variable $x$.
To prove their linear independence we can reduce further the problem to the
case where $V_1^+(x) = \frac 1{d+1}x^{d+1}$. Indeed, suppose that
there exist constants $A_j$ such that 
\be
\mathcal W(z;v_0,...,v_d):= 
\sum_{j=1}^{s_1} A_j\int_{\Gamma_j}\!{\rm d}x\,W_1(x){\rm e}^{xz} \equiv 0\ ,
\label{red2}
\ee
where we have emphasized the dependence on the subleading coefficients of
$V_1^+$ as given in Eqs. (\ref{potent}, \ref{weightx}). Considering
it as a function of the variables $ v_0,...,v_d$
then Eq. (\ref{red2}) implies that 
\be
\frac {\pa^{|\alpha|}}{\pa {\tilde{\underline v}}^\alpha} \mathcal W(z;\tilde {\underline v})\bigg |_{\tilde v_i=v_i}
= 0\ ,\ \ \forall\, \alpha = (\alpha_1,...,\alpha_d)\in \N^d\ ,\
\forall z \in \C.\label{indep2}
\ee
Since $\mathcal W(z;\tilde v_0,....,\tilde v_d)$ is clearly entire
in the variables $\tilde v_i$, Eq. (\ref{indep2})  implies that
actually it does not depend on them. In other words if the $\Xi_i$s
are linearly dependent with constants $A_i$ then also the $\Xi_i$s
where we ``switch off'' the coefficients $v_i$ of the potential are
 linearly dependent with the {\em same} constants $A_i$.\\
Therefore it also does not affect the generality of the problem of
showing linear independence to assume the specific form for $V_1^+$
\be
V_1^+(x)  = \frac 1{d+1}x^{d+1}\ .
\ee
We now analyze the asymptotic behavior, and we need the following
definition  (here given for a $V_1^+$ more general than the one above).
\bd
\label{defSDC}
The steepest descent  contours (SDCs) for integrals of the form  
\be
I_\Gamma(z):= \int_\Gamma\!{\rm d}x\, {\rm e}^{-V_1^+(x)+xz} H(x)\ ,
\ee
with $H(x)$  of polynomial growth at $x=\infty$, are the contours $\gamma_k$
uniquely defined, as $z\to \infty$ within the sector $\mathcal E =
\le\{
\arg(z)\in \le(-\frac{\pi}{2(d+1)},0\ri)\ri\}$, by
\bea
&& \gamma_k:=\le\{x\in \C;\ \Im(V_1^+(x)-xz) =
\Im\le(V_1^+(x_k(z))-zx_k(z)\ri)\  , \Re
(V_1^+(x))\m{\longrightarrow}_{
\shortstack{
\scriptsize $x\to\infty$\\
\scriptsize $x\in \gamma_k$}} +\infty\ .\ri\}\ ,
\eea
where $x_k(z)$ are the $d_1$ branches of the solution to 
\be
{V_1^+}'(x)=z\ ,
\ee
which behave as $ z^{\frac 1{d_1}}$ as
$z\to\infty$ in the sector,  for
the different determinations of the roots of $z$.\\
Their homology class is constant as $x\to\infty$ within the sector.
\ed
With reference to Figure 1, the sector $\mathcal E$ is the narrow 
dark-shaded dual sector of $\mathcal S_L$ (light-shaded).
\bp
\label{asymptotic}
Let $\mathcal E$ be the sector $\arg(z)\in \le(-\frac{\pi}{2(d+1)},0\ri)$ at
 $z=\infty$.
 Then the Laplace-Fourier transforms over the SDCs $\gamma_k$  
\be
F_k(z):= \int_{\gamma_k} \! {\rm d}x \, W_1(x) {\rm e}^{zx} \ ,
k=1,\dots d\
\ee
 have the following asymptotic leading behavior in the sector
$\mathcal E$
\bea
&&F_k(z) =K\sqrt{\frac {2\pi}d}  z^{\frac {2A+1-d}{2d}}
\omega^{k(A-\frac 12)} \exp\le[ \frac d{d+1} z^{\frac {d+1}d}
\omega^k\ri]\le(1+\mathcal O\le(\frac 1 z \ri)\ri)\ ,\ \\
&&A:=\sum_{j=1}^p\lambda_j\ ,\ \omega:= {\rm e}^{\frac{2i\pi}d}\ ,
\eea
where $K\neq 0$ is a constant found in the proof.
\ep
{\bf Proof}.\\
The proof of this asymptotic is an application of the saddle point
method.
Writing $z=|z|{\rm e}^{i\theta}$  with the change
$x=|z|^{1/d}\xi$ we can rewrite the integrals
\bea
&& \int_\Gamma{\rm e}^{-\frac 1{d+1}x^{d+1}+xz}
 \prod_{j=1}^{p} (x-X_j)^{\lambda_j} {\rm e}^{T(x)}
{\rm d}x =\\
&&=|z|^{\frac 1 d} |z|^{\frac Ad} \int_\Gamma 
{\exp}\le[-|z|^{\frac {d+1}d}\le(\frac 1{d+1}\xi^{d+1}
 -\xi{\rm e}^{i\theta} \ri)\ri] \xi^A\prod_{j=1}^p \le(1 - \frac
{X_j}{\xi|z|^{\frac 1d}} \ri)^{\lambda_j} {\rm e}^{T(|z|^{1/d} \xi) }
 {\rm d}\xi \ ,\\
&&T(x):= \exp\le[\frac {M_1(x)}{\prod_{j=1}^p (x-X_j)^{g_j}}\ri] \to
 K\neq 0 \ , |x|\to\infty. 
\eea
Let us set $\lambda:= |z|^{\frac {d+1}d}$ and change integration variable
\be
s=S(\xi):= \frac 1{d+1}\xi^{d+1} 
 -\xi{\rm e}^{i\theta}\ .
\ee
Note that the rescaling of variable leaves the contour $\Gamma$ in
the same ``homology'' class, so that we can take the contour as fixed
in the $\xi$-plane.
The saddle points for this exponential are the roots of 
\be
0= S'(\xi) =\xi^d -{\rm e}^{i\theta} \ ,
\ee  
that is the $d$
roots of ${\rm e}^{i\theta}$.
 The corresponding critical values are
\be
s_{cr}^{(k)}(\theta) := -\frac {d}{d+1} \omega^k{\rm e}^{i\theta\frac{d+1}d} \ ,\ \
\omega :={\rm e}^{2i\pi/d},\ \ k=0,\dots, d-1.\label{critval}
\ee
The map $s=S(\xi)$ is a $d+1$-fold covering of the $s$ plane by the
$\xi$-plane with square-root-type 
branching points at the $s_{cr}^{(k)}(\theta)$.
Moreover each of the $d+1$ sectors (around $\xi=\infty$) for which
$\Re(\xi^{d+1})>0$  is mapped to
the single sector
\be
\mathcal S:= \{s\in \C,\ \ -\frac \pi 2+\epsilon < \arg(s)<\frac \pi
2-\epsilon \}\ .
\ee 
The inverse map $\xi=\xi(s)$ is univalued if we perform the  cuts on the
$s$ plane starting at each $s_{cr}^{(j)}(\theta)$ and going to
$\Re(s)=+\infty$ parallel to the real axis.
Such cuts are distinct for generic values
of $\theta$. We obtain a simply connected domain in the $s$ plane
(see picture).
By their definiton the SDCs $\gamma_j$ corresponds to (the two rims of)
the horizontal cuts  in the
$s$-plane that go from the  critical points $s_{cr}^{(j)}(\theta)$ to
$\Re(s)=+\infty$.\\
The cuts are distinct 
 if $\Im\le({\rm e}^{i\frac {d+1}d\theta
+2ik\frac{\pi}{d}}\ri)\neq \Im\le({\rm e}^{i\frac {d+1}d\theta
+2ij\frac\pi d}\ri) $, for $j\neq k$, that is  away from the Stokes' lines at infinity
\be
{ l}_k = \le\{ \arg(z) = \frac {\pi k}{d+1},\ \ k\in \frac 12 \Z \ri\}.
\ee
Therefore if $z$ approaches infinity along a ray distinct from the
Stokes' lines and within the same sector between them,
 the asymptotic expansion does not change.\\
{\bf Asymptotic evaluation of the steepest descent integrals}\\
The integrals corresponding to the steepest descent path $\gamma_k$  become
\bea
&& |z|^{\frac {A+1}d} \int_{\gamma_k} 
{\rm e}^{-\lambda s} \xi(s)^A g(s,|z|) \frac {d\xi}{ds} 
 {\rm d}s\ ,\\
&&g(s,|z|):= \prod_{j=1}^p \le(1 - \frac
{X_j}{\xi(s) |z|^{\frac 1d}} \ri)^{\lambda_j} {\rm e}^{T(|z|^{1/d}
\xi(s)) },\ \ \lim_{|z|\to \infty} g(s,|z|) = K\neq 0.
\eea
where $\lambda := |z|^{\frac {d+1}d}$.
The  Jacobian of the change of variable has square-root types
singularity at the critical point $s_{cr}^{(k)}$ since 
the singularities (in the sense of singularity theory) 
 of $S(\xi)$ are simple and nondegenerate.\\[3pt]
\centerline{\parbox{8cm}{\epsfxsize=8cm
\epsfysize=4cm
\epsffile{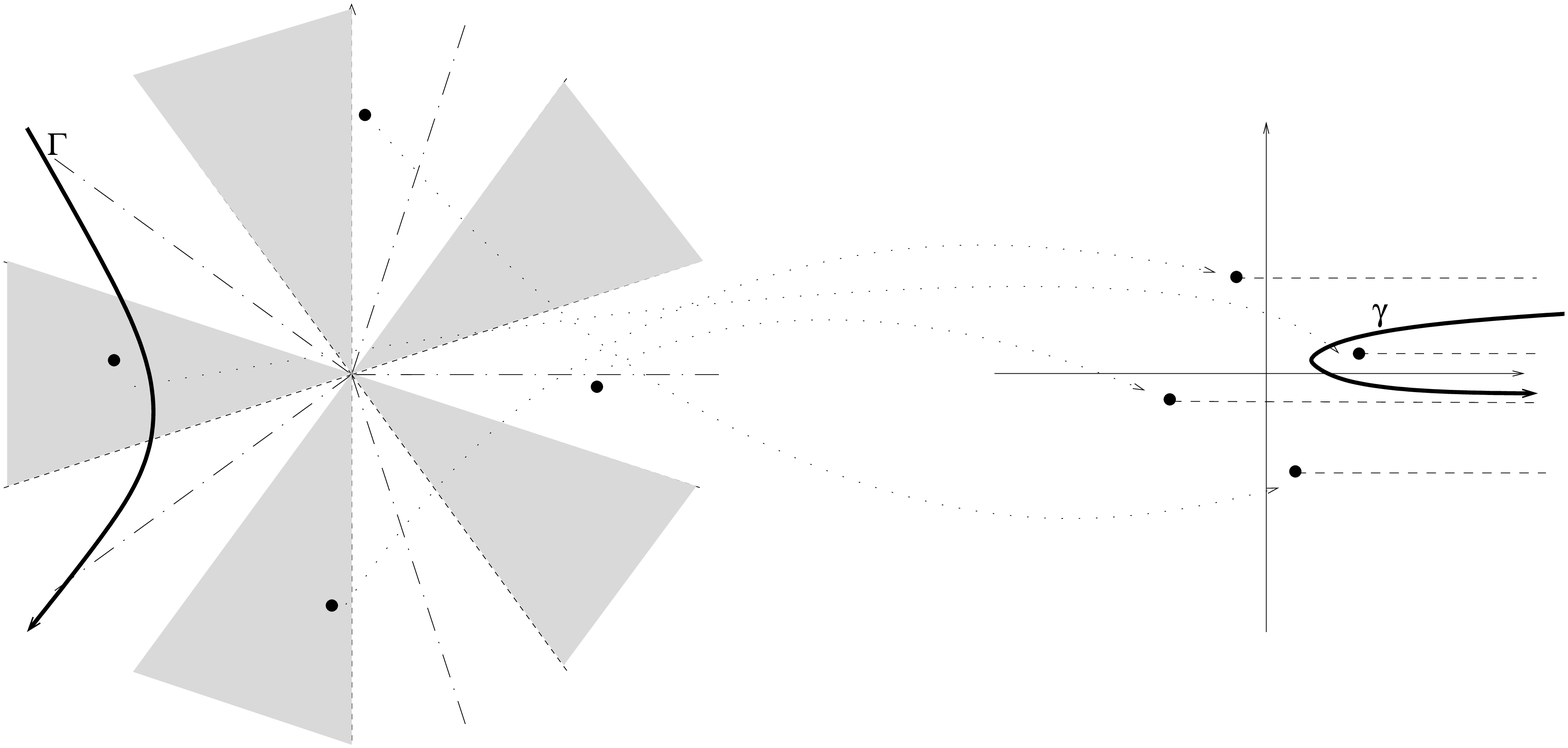}}}
\centerline{Figure 2: The Steepest Descent contours for $d=4$. The
left depicts the $\xi$-plane, the right the $s$-plane.}\\[5pt]
 Then the above integral
becomes, upon developing the Jacobian in Puiseux series,
\bea
&& |z|^{\frac {A+1} d} \int_{\gamma_k} {\rm e}^{-\lambda s}
g(s,|z|)\xi(s)^A\frac {d\xi}{ds}(s)  =\\
&&= 
|z|^{\frac {A+1} d}{\rm e}^{-\lambda s_{cr}}
 \int_{\gamma_k}  {\rm d}s\,{\rm
e}^{-\lambda (s-s_{cr})}\xi(s)^A  g(s,|z|) \frac  {1 } 
{\sqrt{2\frac {d^2s}{d\xi^2}(s_{cr})(s-s_{cr})}}(1+...) = \\
&& \simeq 
K |z|^{\frac {A+1}d} {\rm e}^{i\frac A d \theta}
\omega^{kA} {\rm e}^{-\lambda s_{cr}}
 \le( 2d{\rm e}^{\frac {d-1}d\theta}\omega^k
\ri)^{-\frac 1 2 } 2\int_{\R_+} {\rm e}^{-\lambda t}\frac{{\rm
d}t}{\sqrt{t}} = \\
&&=K |z|^{\frac {A+1} d}\omega^{kA} {\rm e}^{i\frac Ad \theta}
 {\rm e}^{-\lambda s_{cr}} 
\le( 2d{\rm e}^{\frac {d-1}d\theta}\omega^k
\ri)^{-\frac 1 2 } 2\sqrt{\pi}\lambda^{-\frac 12} = \\
&&= K\sqrt{\frac {2\pi}d}  z^{\frac {2A+1-d}{2d}}
\omega^{k(A-\frac 12)} \exp\le[ \frac d{d+1} z^{\frac {d+1}d}
\omega^k\ri]\ .Q.E.D.
\eea
In particular Proposition \ref{asymptotic} shows that the SDC
integrals $F_k$ are linearly independent because their asymptotics
clearly is.\\
Since the SDCs $\gamma_k$  and the contours $\Gamma_k$ span the same
homology, we can always assume that the $\Xi_i$ corresponding to the
closed loops attached to $\infty$ are integrals over the SDC
$\gamma_k$ 
Suppose now that there exist constants $A_i$ such that 
\be
\sum_{j=1}^{s_1} A_i\Xi_i(z) \equiv 0\ .\label{due}
\ee
We split the sum into two parts; the first one contains all contour
integrals 
corresponding to the bounded paths, the paths joining the finite
zeroes $X_i$s to infinity, and loops attached to $X_0=\infty$
approaching $\infty$ within the sector $\mathcal S_L$.
 We denote the subset
of the corresponding  indices by $I_L$. 
Now it is a simple check which we leave to the reader 
that all these integrals are of exponential
type in the sector $\mathcal E$ dual to $\mathcal S_L$\footnote{Saying
that a function is of exponential type in a given
sector means that there exist constants $K$ and $C$ such that the
function is bounded by $|z|^K {\rm e}^{C|z|}$ in that sector.}.\\
The second subset of indices $I_R$ corresponds to the remaining contour integrals
over paths which come from and return to $\infty$ outside the sector
$\mathcal S_L$; a careful counting gives $|I_R| = [d/2]$.
The sum in (\ref{due}) can be accordingly separated in
\bea
\sum_{i\in I_L} A_i \Xi_i(z) = -\sum_{i\in I_R} A_i\Xi_i(z)\ .
\label{split}
\eea
We want to conclude that the two sides of Eq. (\ref{split}) must
vanish separately.
Indeed we have remarked above that 
the LHS in (\ref{split}) is of exponential type in the
sector $\mathcal E$.\\
On the other hand we now prove that 
the RHS  {\em cannot} be of exponential type unless each
of  $A_i,\ i\in I_R$ vanishes.
From Prop. \ref{asymptotic} we deduce that among the SDC
integrals there are precisely $[d/2]$ that have a dominant exponential
behavior of the type $
\exp\le( \frac d{d+1} z^{\frac {d+1}d}
\omega^{k}\ri)$ with $\Re( z^{\frac {d+1}d}
\omega^{k})>0$ in the sector $\mathcal E$, whic is {\em not} of
exponential type; since the SDC's can be
obtained by suitable linear combinations with integer coefficients of
the chosen contours then the $[d/2]$ functions $\Xi_i,\ i\in
I_R$ must span the same space as the dominant $[d/2]$ linearly
independent SDC's in the sector $\mathcal E$, {\em modulo} the span of
$\Xi_i,\ i\in I_L$. In formulae
\be
\Z\{F_k:\ F_k \hbox{ dominant in }\mathcal E\}\  \simeq\ 
 \Z\{\Xi_i,\
\forall i\}{\rm mod}\, \Z\{\Xi_i,\ i\in  I_L\}\  =\ 
\Z\{\Xi_i,\ i\in I_R\}\ .
\ee 
Since no nontrivial linear combination of the  $[d/2]$ dominant SDC
integrals $F_k$'s in $\mathcal E$ can be of exponential type, the  only possibility for the RHS of Eq.
(\ref{split}) to be of exponential type in the sector $\mathcal E$ 
is that $$A_i=0,\ \forall i\in I_R\ .$$
Let us now focus on the terms in the LHS of Eq. (\ref{split}).
We must now  prove that also $A_i=0,\ i\in I_L$. We can now follow
\cite{miller} without hurdles. 
We sketch the main steps below for the sake of completeness.\\
We need to prove that 
\be
Q(z):= \sum_{i\in I_L} A_i\int_{\Gamma_i}{\rm d}x\, W_1(x)\, {\rm e}^{xz}\equiv
0\  \Leftrightarrow\ A_i=0 \ \forall i\in I_L.\label{master}
\ee
Let $a$ be a point within the sector $\mathcal E$ and far enough from
the origin so as to leave all contours $\Gamma_i,\ i\in I_L$ to the
left\footnote{More precisely in the half plane to the left of the
perpendicular to the bi-secant of the sector $\mathcal E$}.
 Let us choose a contour $\mathcal C$ starting at $z$ and going
to infinity in the sector to $\mathcal E$.
Then we integrate  $Q(\zeta) {\rm e}^{-a\zeta}$ 
along $\mathcal C$. Since ${\rm e}^{\zeta (x-a)} W_1(x)$ is
jointly absolutely integrable with respect to the arc-length on each
of the $\Gamma_i,\ i\in I_L$ and $\mathcal C$, we may interchange the order of
integration to obtain 
\be
\sum_{i\in I_L}A_i
 \int_{\Gamma_i}\frac 1{x-a}{\rm e}^{z(x-a)}  W_1(x) \equiv 0. 
\ee
Repeating this $r-1$ times and then setting $z=0$ at the end,  we
obtain 
\be
\sum_i A_i\int_{\Gamma_i} \big(x-a\big)^{-r} W_1(x){\rm d}x \equiv
0,\ \forall r\in \N.\label{master2}
\ee
Let us define 
\be
\tilde v(x):=W_1(x) (x-a)^2
\ee
so that Eq. (\ref{master2}) is turned into 
\be
\sum_i A_i \int_{\Gamma_i} \big(x-a\big)^{-r} \tilde v(x) \frac{{\rm
d}x}{(x-a)^2} \equiv
0,\ \forall r\in \N.\label{master22}
\ee
Let us perform the change of variable $\omega = \frac 1{x-a}$ (a homographic
transformation). We denote by $\gamma_i$ the images of the contours
$\Gamma_i$ and by $f(\omega)$ the function $\tilde v(x(\omega))$. \\ 
Eq. (\ref{master2}) (or equivalently Eq. (\ref{master22})) now becomes 
\be
\sum_{i\in I_L} A_i\int_{\gamma_i} \!{\rm d}\omega f(\omega) P(\omega) = 0 \
,\ \ \forall P\in \C[\omega]\ .\label{master3}
\ee 
Note that in the variable $\omega$ all contours are in the finite
region of the $\omega$-plane and the contours look like the ones in
Figure 3 (the missing loops attached to $0 = \omega(X_0) =
\omega(\infty)$ were the contours indexed by $I_R$).\par
\centerline{\parbox{10cm}{
\epsfxsize=8cm
\epsfysize=6cm 
\epsffile{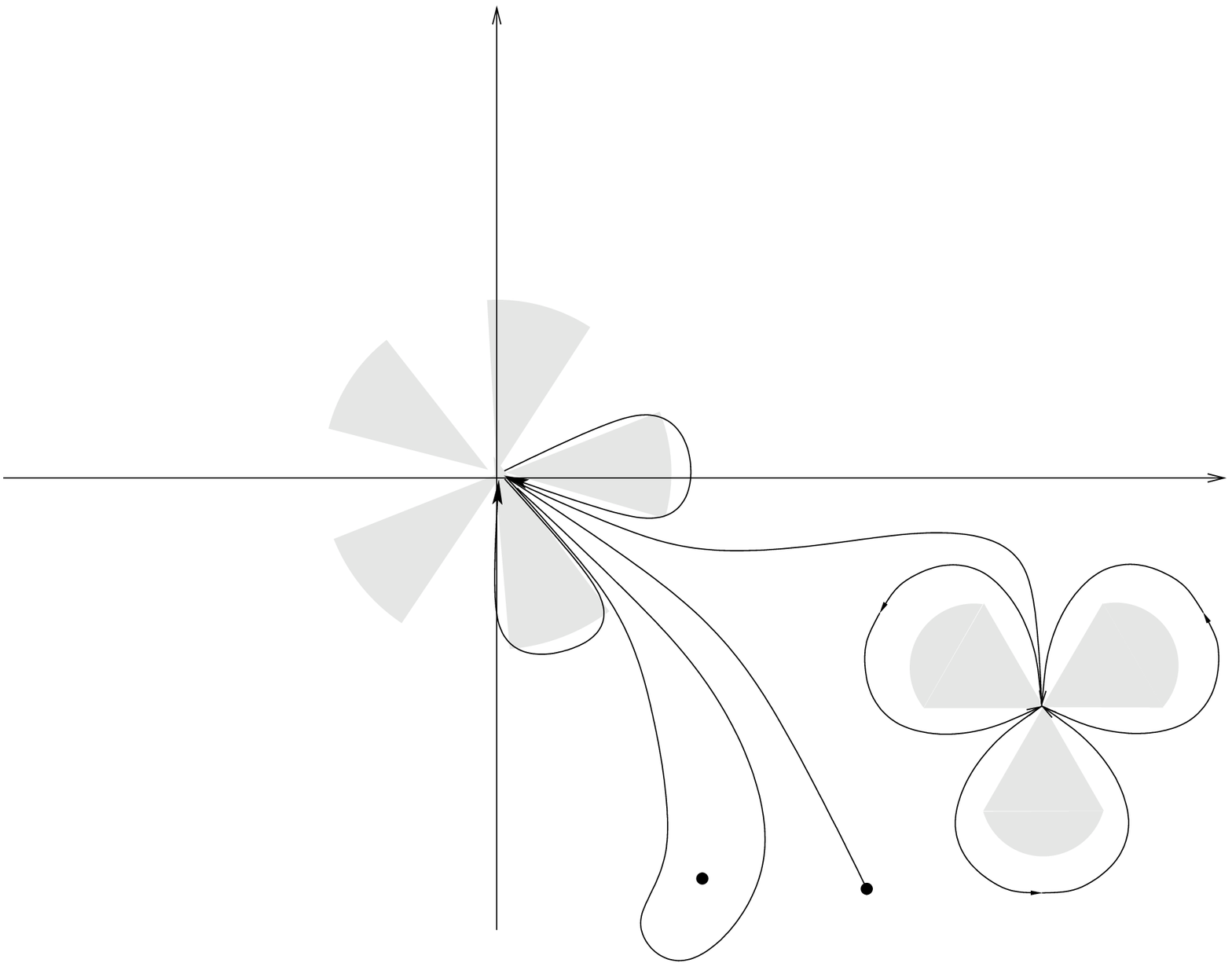}
}}
\centerline{Figure 3: The contours $\gamma_i,\ i\in I_L$ in the
$\omega$ plane.}
We denote by $E$ the closed and bounded set in the $\omega$ plane
constituted by all contours $\gamma_i$, $i\in I_L$  and the interiors of the closed
loops. This set $E$ satisfies the requirements of Lemma (\ref{margelyan}).
Moreover the contours $\gamma_i$ have all the Property
($\wp$) with respect to $f(\omega)$. \par
We now start proving that the $A_i$s vanish.\\
First consider a contour $\gamma_i$ without interior points
(i.e. those segments which join two different $X_i$s). Let $\omega(t)$
be a parametric representation where $t\in [0,L]$ is the arc length parameter
parameter so that $\omega'(t)$ is continuous and nonvanishing on
$[0,L]$. Therefore it follows that the function 
\be
\chi_i(\omega):= \le\{\begin{array}{ll}
\frac {\overline{f(\omega)}} {\omega'(t)}, &\omega \in \gamma_i\cr
0, & \omega\in E\setminus \gamma_i
\end{array}\ri.
\ee
is continuous on $E$ and analytic in the interior points of $E$. Hence
there exists a sequence of polynomials $P_n(\omega)$ converging
uniformly to $\chi_i(\omega)$ on $E$ (by Lemma
\ref{margelyan}). Plugging into Eq. (\ref{master3}) and passing to the
limit we obtain
\be
A_i\int_0^L {\rm d}t\, |f(\omega(t))|^2 = 0\ ,
\ee
which implies that $A_i$ vanishes.\par
Let us now consider a closed loop, say $\gamma_l$. 
Let $T(\omega)$ be any polynomial vanishing at $\omega_0\in \gamma_l$
where $\omega_0$ is the image of the (unique) zero of $B_1(x)$ on the
contour $\Gamma_l$. Then we define 
\be
\Phi_l(\omega):= \le\{
\begin{array}{ll}
T(\omega), & \omega\in \gamma_l \hbox{ and its interior}\\
0, & \omega\in E\setminus\{\gamma_l \hbox{ and its interior}\}
\end{array}
\ri.
\ee
Again, $\phi_l(\omega)$ satisfies the requirement of Lemma
(\ref{margelyan}) and hence can be approximated uniformly by a
sequence of polynomials. Passing the limit under the integral we then
obtain 
\be
A_l \int_{\gamma_l}\!{\rm d}\omega f (\omega) T(\omega) = 0\ ,
\ \  \forall\,
T\in (\omega-\omega_0)\C[\omega]\ .
\ee 
We then use Theorem \ref{millsha} to conclude that $f$ should be
bounded inside $\gamma_l$. But this is a contradiction because $f(\omega)$
has the Property ($\wp$) w.r.t. $\gamma_l$ since  $\tilde v(x) =
W_1(x)(x-a)^2$ had the same Property w.r.t. the
closed contour $\Gamma_l$. This is a contradiction unless the
 $A_l$ vanishes.\par
Therefore we have proven that all the $A_i$ must vanish, i.e. the
$\Xi_i(z)$ are linearly independent.\par
Repeating for the $\Psi_j(w)$ we conclude the proof of Theorem
\ref{main}.
\section{Conclusion}  
We  make a few remarks on the cases
we have not considered, i.e. when $\deg(A_i)\leq \deg(B_i)$ for one or
both $i=1,2$. Indeed (up to some care in the definition of the
contours for reasons of convergence) one can easily define {\em some}
solutions of Eqs. (\ref{oversys}) in the form of double
Laplace--Fourier integrals and also prove their linear
independence. More complicated is to produce the analog of
Prop. \ref{dimension}, that is to have an a-priori knowledge of the
dimension of the space of solutions to Eqs. (\ref{oversys}):
 the result (which we do not prove here) is that there are
$M=s_1s_2+1$ solutions. The moment recurrences
(\ref{momrec1},\ref{momrec2}) say then that the bifunctionals are
actually $M-1$ in Case AB or $M-2$ in Case AA. That is one has to give
a criterion to select amongst the solutions to
Eq. (\ref{oversys}) the ones which are analytic at $w=0=z$. We will
return on this point in a future publication.\\
Suffices here to say that a 
similar problem occurs for the semi-classical moment
functionals $\L:\C[x]\to \C$. As we have illustrated in the
introduction the generating function satisfies Eq. \ref{genfu},
but in general not all solutions are analytic at $z=0$
and hence do not define any moment functional. This can be
understood by looking at the recurrence relations satisfied by the
moments:
\be
n\sum_{j=0}^{d}\beta(j)\mu_{n+j-1} = \sum_{j=0}^{k}\alpha(j)\mu_{n+j}\ ,
\ee 
where $d=\deg(B)> \deg(A)+1 = k+1$. In this case the resulting $d$-terms
recurrence relation has actually only $d-1$ solutions because, for
$n=0$ the above equation gives a {\em constraint} on the initial
conditions\footnote{When $\deg(A)+1=\deg(B)=d$ then {\em generically}
there are $d-1$ solutions, except in some cases when
$\exists n\ s.t.\ \a(d-1) = n\b(d) $. See \cite{marce1} for more details.}
\be
0=  \sum_{j=0}^{k}\alpha(j)\mu_{j}\ .
\ee
This should be regarded as the requirement that the solution of
Eq. (\ref{genfu}) be analytic at $z=0$.\\
Now, in the bilinear case we have the additional problem that the
recurrence relations for the bi-moments are overdetermined and hence
the corresponding constraint on the initial conditions must be shown
to be compatible as well. We postpone the more detailed discussion of
this problem to a future publication.\par\vskip 5pt
{\large \bf Acknowledgments}\\
The author wishes to thank Prof. B. Eynard and Prof. J. Harnad for stimulating
discussion, and Prof. H. S. Shapiro for helpful hints in amending
the proof of linear independence.  

\end{document}